\documentclass{elsart}
\usepackage{amssymb}
\usepackage{graphicx}

\begin{document}

\begin{frontmatter}

\title{Global bifurcations of multiple limit cycles in~the~FitzHugh--Nagumo~system\thanksref{label1}}

\thanks[label1]{The author is very grateful to the Max Planck Institute for Mathematics (Bonn)
for hospitality and support during his stay in March\,--\,April 2011 at this Institute.}

\author{Valery A. Gaiko}
\ead{valery.gaiko@yahoo.com}
\address{United~Institute~of~Informatics~Problems,~National~Academy~of~Sciences~of~Belarus,
Leonid~Beda~Str.\,6-4,~Minsk~220040,~Belarus}

\begin{abstract}
In this paper, we complete the global qualitative analysis of the well-known Fitz\-Hugh--Nagumo neuronal model.
In particular, studying global limit cycle bifurcations and applying the Wintner--Perko termination principle
for multiple limit cycles, we prove that the corresponding dynamical system has at most two limit cycles.
\end{abstract}

\end{frontmatter}

\section{Introduction}

We consider the well-known FitzHugh--Nagumo model in the form
    $$
    \begin{array}{l}
\dot{V}=I-W-aV+(a+1)V^2-V^3,
    \\
\dot{W}=\varepsilon(V-\delta\,W),
    \end{array}
    \eqno(1.1)
    $$
where $V$ is the membrane potential, $W$ is a re\-co\-very variable, and $I$ is the magnitude of stimulus current,
which is a two-dimensional simplification of the classical Hodgkin--Huxley model of the spike dynamics in
a biological neuron \cite{fh}, \cite{Izhikevich}, \cite{nay}, \cite{rzh}, \cite{RGG}.
This system was suggested by FitzHugh (1961) \cite{fh}, who called it ``Bonhoeffer\,--\,van der Pol model'', and the equivalent
circuit was constructed by Nagumo et al. (1962) \cite{nay}. The motivation for the FitzHugh--Nagumo model was to isolate conceptually
the essentially mathematical pro\-per\-ties of excitation and propagation from the electrochemical properties of sodium and potassium ion flow. The model consists of a voltage-like variable having cubic nonlinearity that allows regenerative self-excitation via a positive
feedback, and a recovery variable having a linear dynamics that provides a slower negative feedback. While the Hodgkin--Huxley
model is more realistic and biophysically sound, only projections of its four-dimensional phase trajectories can be observed.
The simplicity of the FitzHugh--Nagumo model permits the entire solution to be viewed at once. This allows a geometrical
explanation of important biological phenomena related to neuronal excitability and spike-generating mechanism~\cite{Izhikevich}.
    \par
The phase portrait of the FitzHugh--Nagumo model (1.1) depicts the $V$-null\-cline, which is the $N$-shaped curve obtained
from the condition $\dot{V}=0,$ and the $W$-null\-cline, which is a straight line obtained from the condition $\dot{W}=0.$
The intersection of null\-clines is an equilibrium (a singular point) of the system (1.1), which may be unstable if it is
on the middle branch of the $V$-null\-cline, i.\,e., when $I$ is strong enough. In this case, the model exhibits periodic
(tonic spiking) activity.
    \par
The FitzHugh--Nagumo model explained the absence of all-or-none spikes in the Hodgkin--Huxley model in response to stimuli, i.\,e.,
pulses of the injected current $I.$ Weak stimuli (small pulses of $I)$ result in small-amplitude trajectories that correspond
to subthreshold responses; stronger stimuli result in intermediate-amplitude trajectories that correspond to partial-amplitude spikes;
and strong stimuli result in large-amplitude trajectories that correspond to suprathreshold response --- firing a spike.
    \par
Similarly to the Hodgkin--Huxley model, the FitzHugh\,--\,Nagumo model does not have a well-defined firing threshold in the absence
of a saddle equilibrium. This feature is the consequence of the absence of all-or-none responses. The apparent illusion of threshold dynamics and all-or-none responses in both models is due to the existence of the ``quasi-threshold'', which is a canard trajectory
that follows the unstable (middle) branch of the $N$-shaped $V$-nullcline.
    \par
The FitzHugh--Nagumo model explains the excitation block phenomenon, i.\,e., the cessation of repetitive spiking as the amplitude of
the stimulus current increases. When $I$ is weak or zero, the equilibrium (intersection of nullclines) is on the left (stable) branch of $V$-nullcline, and the model is resting. Increasing $I$  shifts the nullcline upward and the equilibrium slides onto the middle (unstable) branch of the nullcline. The model exhibits periodic spiking activity in this case. Increasing the stimulus further shifts the equilibrium to the right (stable) branch of the $N$-shaped nullcline, and the oscillations are blocked (by excitation). The precise mathematical mechanism involves appearance and disappearance of a limit cycle attractor, and it is reviewed in detail by Izhikevich (2007) \cite{Izhikevich}.
    \par
This model explained also the phenomenon of post-inhibitory (rebound) spikes, called anodal break excitation at that time. As the stimulus $I$ becomes negative (hyperpolarization), the resting state shifts to the left. As the system is released from hyperpolarization (anodal break), the trajectory starts from a point far below the resting state (outside the quasi-threshold), makes a large-amplitude excursion, i.\,e., fires a transient spike, and then returns to the resting state.
    \par
The FitzHugh--Nagumo model explained the dynamical mechanism of spike accommodation in Hodgkin--Huxley-type models. When stimulation strength  increases slowly, the neuron remains quiescent. The resting equilibrium of (1.1) shifts slowly to the right, and the state of the system follows it smoothly without firing spikes. In contrast, when the stimulation is increased abruptly, even by a smaller amount, the trajectory could not go directly to the new resting state, but fires a transient spike; see figure. Geometrically, this phenomenon
is similar to the post-inhibitory (rebound) response.
    \par
The FitzHugh--Nagumo equations became a favorite model for reaction-diffusion systems
    $$
    \begin{array}{l}
\dot{V}=I-W-aV+(a+1)V^2-V^3+V_{xx},
    \\
\dot{W}=\varepsilon(V-\delta\,W),
    \end{array}
    \eqno(1.2)
    $$
which simulate propagation of waves in excitable media, such as heart tissue or nerve fiber. Here, the diffusion term $V_{xx}$ is the second derivative with respect to the spatial variable $x.$ Its success is mostly due to the fact that the model is analytically tractable, and hence it allows derivation of many important properties of traveling pulses without resort to computer simulations.
	\par	
Without loss of generality, the system (1.1) can be written in the canonical form
    $$
    \begin{array}{l}
\dot{x}=(\gamma\,\delta-1)\,y+(\gamma-a)\,x+b\,x^2-c\,x^3\equiv P(x,y),
    \\
\dot{y}=x-\delta\,y\equiv Q(x,y).
    \end{array}
    \eqno(1.3)
    $$
Such a system was studied earlier, e.\,g., in \cite{rzh}. However, its qualitative analysis
was incomplete, since the global bifurcations of multiple limit cycles could not be studied properly by means
of the methods and techniques which were used earlier in the qualitative theory of dynamical systems.
Applying new bifurcation methods and geometric approaches developed in \cite{bg}, \cite{Gaiko}--\cite{gai6},
we complete the qualitative analysis of the FitzHugh--Nagumo model and prove, in particular, that the corresponding
dynamical system (1.3) has at most two limit cycles. In Sections~2\,--\,4 of this paper, we recall basic facts and results
from the global bifurcation theory of polynomial dynamical systems and its applications. These results, together with the
methods of \cite{bog}, \cite{bg}, \cite{Gaiko}--\cite{gai6}, are used in Sections~5,~6 for the study of singular point and
limit cycle bifurcations of the system~(1.3).

\section{Preliminaries}

In this paper, geometric aspects of Bifurcation and Catastrophe Theories are used and de\-ve\-loped \cite{Gaiko}, \cite{Kuznetsov}, \cite{Perko}. First of all, the two-isocline method which was developed by Erugin is used, see \cite{Gaiko}. An isocline portrait is
the most na\-tu\-ral construction for a polynomial equation. It is sufficient to have only two nullclines (or isoclines of zero and infinity in our terminology) to obtain principal information on the original polynomial system, because these two isoclines are
right-hand sides of the system. Geometric properties of isoclines (conics, cubics, quartics, etc.) are well-known, and all isocline portraits can be easily constructed. By means of them, all topologically different qualitative pictures of integral curves to within
a number of limit cycles and distinguishing center and focus can be obtained. Thus, it is possible to carry out a rough topological classification of the phase portraits for the polynomial dynamical systems. It is the first application of Erugin's method. After studying contact and rotation properties of the isoclines, the simplest (canonical) systems containing limit cycles can be also constructed. Two groups of parameters can be distinguished in such systems: static and dynamic. Static parameters determine the behavior of phase trajectories in principle, since they control the number, position, and character of singular points in a finite part of the plane (finite singularities). The parameters from the first group determine also a possible behavior of separatrices and singular points at infinity (infinite singularities) under variation of the parameters from the second group. The dynamic parameters are field rotation parameters, see \cite{BL}, \cite{Gaiko}, \cite{Perko}. They do not change the number, position and index of the finite singularities, but only involve the vector field in a directional rotation. The rotation parameters allow to control the infinite singularities, the behavior of limit cycles and separatrices. The cyclicity of singular points and separatrix cycles, the behavior of semi-stable and other multiple limit cycles are controlled by these parameters as well. Therefore, by means of the rotation parameters, it is possible to control all limit cycle bifurcations and to solve the most complicated problems of the qualitative theory of dynamical systems.
    \par
In \cite{Gaiko}, \cite{gai1}, \cite{gai2}, \cite{gai3}, \cite{gai6}, some complete results on quadratic systems have been presented.
In particular, it has been proved that for quadratic systems four is really the maximum number of limit cycles and $(3\!:\!1),$ i.\,e., three limit cycles around one focus and the only limit cycle around another focus, is their only possible distribution (this is a solution of Hilbert's Sixteenth Problem in the quadratic case of polynomial dynamical systems). In \cite{gvh}, some preliminary results on generalizing new ideas and methods of \cite{Gaiko} to cubic dynamical systems have already been established. In particular, a canonical cubic system of Kukles type has been constructed and the global qualitative analysis of its special case corresponding to a generalized Li\'{e}nard equation has been carried out. It has been proved also that the foci of such a Li\'{e}nard system can be at most of second order and that such system can have at most three limit cycles on the whole phase plane. Moreover, unlike all previous works on the Kukles-type systems, global bifurcations of limit and separatrix cycles using arbitrary (including as large as possible) field rotation parameters of the canonical system have been studied in \cite{gvh}. As a result, the classification of all possible types of separatrix cycles for the generalized Li\'{e}nard system has been obtained and all possible distributions of its limit cycles have been found. In \cite{gai4}, \cite{gai5}, a solution of Smale's Thirteenth Problem proving that the Li\'{e}nard system with a polynomial of degree $2k+1$ can have at most $k$ limit cycles has been presented. In \cite{bg}, we have completed the global qualitative analysis of a quartic ecological model. All of these methods and results can be applied to the global qualitative analysis of the Fitz\-Hugh--Nagumo
neuronal model as~well.
    \par
In \cite{bog}, we have already carried out the global qualitative analysis of a polynomial dynamical system as a learning model of
neural networks \cite{bj}, \cite{oja}. Learning models are algorithms, implementable as neural networks, that aim to mimic an adaptive procedure. A~neural network is a device consisting on interconnected processing units, designated neurons. An input presented to the network is translated as a numerical assignment to each neuron. This will create a sequence of internal adjustments leading to a learning process. An input vector, denoted by $\xi,$ represents an $n$-dimensional random vector with independent components. This means that the joint probability distribution function is the product of $n$ density functions. The output value, denoted by~$V,$ is the outcome of the network's action on $\xi$ and is given by $\sum_{j=1}^n\omega_j\xi_j,$ where $\omega_j$ is the connecting weight for the synapse attached to the input neuron $\!j.$ Since new synapses may be created under a constant error rate, $E,$ a synaptic strength may capture nearby
activity. This is done by the creation of temporary synapses from the closest neurons to the output one. The synaptic rate of change
is given by
    $$
    \dot{\omega}_i
    =V((1-E)\xi_i+(E/2)(\xi_{i+1}+\xi_{i-1})-V\omega_i),
    $$
for $i=2,\ldots,\,n-1,$ or
    $$
    \dot{\omega}_i
    =V((1-E)\xi_i+(E/2)\xi_{i\pm 1}-V\omega_i),
    $$
for $i=1$ or $n,$ respectively. Substituting the value of $V$ in the expression of $\dot{\omega}_i,$ we obtain
    $$
    \dot{\omega}_i=
    \left\{
    \begin{array}{l}
    (1\!-\!E)\sum_{j=1}^n\omega_j \xi_j \xi_i+(E/2)\sum_{j=1}^n\omega_j\xi_j
    (\xi_{i-1}+\xi_{i+1})\!-\sum_{j,\,k}\omega_j\omega_k\xi_j\xi_k\omega_i,\\
    i\neq1~\mbox{and}~n,
    \\[2mm]
    (1\!-\!E)\sum_{j=1}^n \omega_j \xi_j \xi_i+(E/2)\sum_{j=1}^n
    \omega_j\xi_j\xi_{i\pm1}\!-\sum_{j,\,k}\omega_j\omega_k\xi_j\xi_k\omega_i,\\
    i=1~\mbox{or}~n,
    \end{array}
    \right.
    $$
what can be reduced to the equation
    $$
    \frac{d{\omega}}{dt}=TC\omega-(\omega,C\omega)\omega,
    \eqno(2.1)
    $$
with $C=[\langle\xi_i\xi_j\rangle]_{ij}=\xi^t\xi,$ \
$\xi=\{\xi_1,\xi_2,\,\ldots,\,\xi_n\},$ a correlation matrix of
expected values, and $T,$ a~tridiagonal substochastic matrix given
by $t_{ij}=0$ if $|i-j|>1,$ $t_{ij}=E/2$ if $|i-j|=1$, and
$t_{ii}=1-E,$ for all $i$ and $j=1,\,\ldots,n$~\cite{bj}.
    \par
Desirable initial conditions are those with trajectories that
converge to some equilibrium (singular) point of~(2.1). This will
assume a natural weight assignment as a result of the learning
process. Knowledge on the stability of equilibria provides
information on the robustness of the learning process. Existence
of cycles might represent a different kind of learning where a
whole continuum of connecting weight vectors emerges instead of
just a single vector. In \cite{bog}, we have restricted our attention
to two dimensions.
    \par
For two input neurons, (2.1) can be written as a cubic dynamical system
    $$
    \begin{array}{l}
 \dot{x}=((1\!-\!\varepsilon)a\!+\!(\varepsilon/2)b)x\!+\!((1\!-\!\varepsilon)
 b\!+\!(\varepsilon/2)c)y-x(ax^2\!+\!2bxy\!+\!cy^2),
    \\[2mm]
 \dot{y}=((\varepsilon/2)a\!+\!(1\!-\!\varepsilon)b)x\!+\!((\varepsilon/2)b\!+
 \!(1\!-\!\varepsilon)c)y-y(ax^2\!+\!2bxy\!+\!cy^2),
    \end{array}
    \eqno(2.2)
    $$
where the parameters $\varepsilon$ and $a,$ $b,$ $c$ represent,
respectively, the probability of synaptic formation and the weight
strengths for the synapses attached to the input neurons~\cite{bog}.
Thus, we have got a four-parameter planar dynamical system for investigation.
Applying techniques based both on classical Poincar\'{e} and Dulac methods
and also on some methods developed in~\cite{Gaiko}, we have studied the
global bifurcations of singular points and limit cycles of the cubic
system~(2.2), a learning model of planar neural networks~\cite{bog}.
    \par
Some of these techniques can be extended to higher-dimensional
dynamical systems \cite{Gaiko}, \cite{Kuznetsov}, \cite{Perko}.
So, for the global analysis of limit cycle bifurcations (in particular,
for solving the uniqueness problem) we have used the Perko planar termination
principle stating that the maximal one-parameter family of multiple limit
cycles terminates either at a singular point, which is typically
of the same multiplicity, or on a separatrix cycle, which is also
typically of the same multiplicity \cite{Perko}. This principle is
a consequence of the Wintner principle of natural termination, which
was stated for higher-dimensional dynamical systems (see \cite{Gaiko},
\cite{Perko}), where one-parameter families of periodic orbits of the
restricted three-body problem are studied and Puiseux series are used
to show that in the analytic case any one-parameter family of periodic
orbits can be uniquely continued through any bifurcation except a
period-doubling bifurcation. Thus, the Wintner--Perko termination principle
and the method developed in \cite{bog}, \cite{bg}, \cite{Gaiko}--\cite{gai6}
can be applied to the further global qualitative analysis of neural dynamical
systems.

\section{Basic facts on limit cycles}

Consider a polynomial dynamical system in the vector form
    $$
    \mbox{\boldmath$\dot{x}$}=\mbox{\boldmath$f$}
    (\mbox{\boldmath$x$},\mbox{\boldmath$\mu$)},
    \eqno(3.1)
    $$
where $\mbox{\boldmath$x$}\in\textbf{R}^2;$ \
$\mbox{\boldmath$\mu$}\in\textbf{R}^n;$ \
$\mbox{\boldmath$f$}\in\textbf{R}^2$ \ $(\,\mbox{\boldmath$f$}$ is
a polynomial vector function).
    \par
Let us recall some basic facts concerning limit cycles of (3.1). But first
let us state two fundamental theorems from the theory of ana\-ly\-tic
functions~\cite{Gaiko}.
    \par
    \medskip
    \textbf{Theorem 3.1 (Weierstrass Preparation Theorem).}
    \emph{Let $F(w,z)$ be an analytic in the neighborhood of the point
$(0,0)$ function satisfying the following conditions}
    $$
    F(0,0)=0, \: \frac{\partial F(0,0)}{\partial w}=0, \:
    \ldots, \: \frac{\partial^{k-1}F(0,0)}{\partial^{k-1}w}=0; \quad
    \frac{\partial^{k}F(0,0)}{\partial^{k}w}\neq0.
    $$
    \par
    \emph{Then in some neighborhood $|w|<\varepsilon,$ $|z|<\delta$ of
the points $(0,0)$ the function $F(w,z)$ can be represented as}
    \vspace{-2mm}
    $$
    F(w,z)=(w^{k}+A_{1}(z)w^{k-1}+\ldots+A_{k-1}(z)w+A_{k}(z))\Phi(w,z),
    $$
\emph{where $\Phi(w,z)$ is an analytic function not equal to zero
in the chosen neighborhood and $A_{1}(z),\ldots,A_{k}(z)$ are
analytic functions for $|z|<\delta.$}
    \par
    \medskip
From this theorem it follows that the equation $F(w,z)=0$ in a sufficiently
small neighborhood of the point $(0,0)$ is equivalent to the equation
    $$
    w^{k}+A_{1}(z)w^{k-1}+\ldots+A_{k-1}(z)w+A_{k}(z)=0,
    $$
which left-hand side is a polynomial with respect to~$w.$ Thus,
the Weierstrass preparation theorem reduces the local study of the
general case of implicit function $w(z),$ defined by the
equation $F(w,z)=0,$ to the case of implicit function,
defined by the algebraic equation with respect to~$w.$
    \par
    \medskip
    \textbf{Theorem 3.2 (Implicit Function Theorem).}
    \emph{Let $F(w,z)$ be an analytic function in the neighborhood of
the point $(0,0)$ and $F(0,0)\!=\!0,$ $F'_{w}(0,0)\!\neq\!0.$}
    \par
    \emph{Then there exist $\delta>0$ and $\varepsilon>0$ such that
for any~$z$ satisfying the condition $|z|<\delta$ the equation
$F(w,z)=0$ has the only solution $w=f(z)$ sa\-tisfying the
condition $|f(z)|<\varepsilon.$ The func\-tion $f(z)$ is
expanded into the series on positive integer powers of~$z$ which
converges for $|z|<\delta,$ i.\,e., it is a single-valued analytic
function of~$z$ which vanishes at $z=0.$}
    \medskip
    \par
Assume that the system (3.1) has a limit cycle
    $$
    L_0:\mbox{\boldmath$x$}=\mbox{\boldmath$\varphi$}_0(t)
    $$
of minimal period $T_0$ at some parameter value
$\mbox{\boldmath$\mu$}=\mbox{\boldmath$\mu$}_0\in\textbf{R}^n$
(Fig.~1).
    \par
    \begin{figure}[htb]
\begin{center}
    \includegraphics[width=138.5mm]{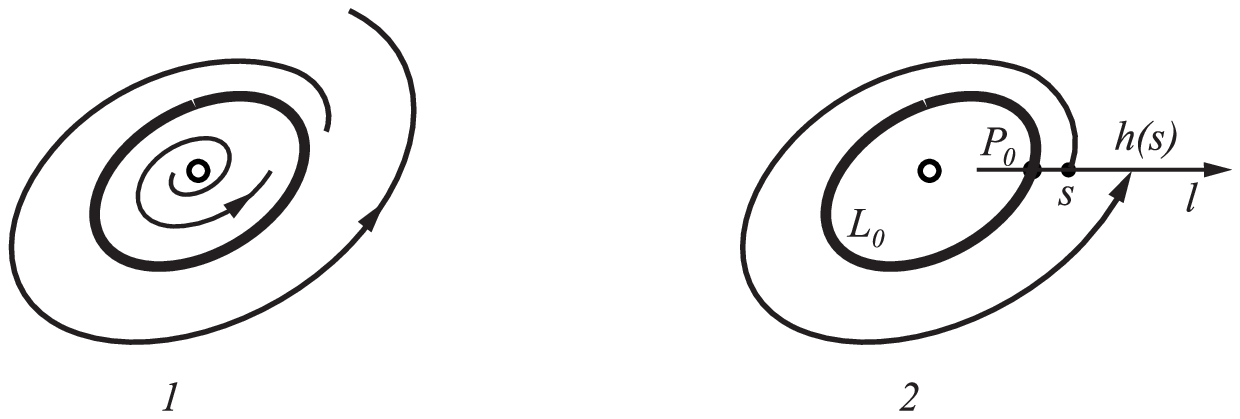}
    \vspace{-4mm}
    \par
    {\small FIG.~1. The Poincar\'{e} return map
    in the neighborhood of a multiple limit cycle.}
\end{center}
    \end{figure}
    \medskip
    \par
Let~$l$ be the straight line normal to $L_0$ at the
point $\mbox{\boldmath$p$}_0=\mbox{\boldmath$\varphi$}_0(0)$ and
$s$ be the coordinate along $l$ with $s$ positive exterior of
$L_0.$ It then follows from the implicit function theorem that
there is a $\delta>0$ such that the Poincar\'{e} map
$h(s,\mbox{\boldmath$\mu$})$ is defined and analytic
for $|s|<\delta$ and
$\|\mbox{\boldmath$\mu$}-\mbox{\boldmath$\mu$}_0\|<\delta.$
Besides, the displacement function for the system
(3.1) along the normal line~$l$ to~$L_0$
is defined as the function
    $$
    d(s,\mbox{\boldmath$\mu$})=h(s,\mbox{\boldmath$\mu$})-s.
    $$
    \par
In terms of the displacement function, a multiple limit cycle can
be defined as follows \cite{Gaiko}.
    \par
    \medskip
    \textbf{Definition 3.1.}
A limit cycle $L_0$ of (3.1) is a \emph{multiple limit cycle} iff
$d(0,\mbox{\boldmath$\mu$}_0)\!=\!d_r(0,\mbox{\boldmath$\mu$}_0)\!=\!0$
and it is a \emph{simple limit cycle} (or hyperbolic limit cycle)
if it is not a multiple limit cycle; furthermore, $L_0$ is a limit
cycle of multiplicity~$m$ iff
    $$
    d(0,\mbox{\boldmath$\mu$}_0)=d_r(0,\mbox{\boldmath$\mu$}_0)=\ldots
    =d_r^{(m-1)}(0,\mbox{\boldmath$\mu$}_0)=0, \quad
    d_r^{(m)}(0,\mbox{\boldmath$\mu$}_0)\neq 0.
    $$
    \par
Note that the multiplicity of $L_0$ is independent of the point
$\mbox{\boldmath$p$}_0\in L_0$ through which we take the normal
line~$l.$
    \par
Let us write down also the following formulas which have already
become classical ones and determine the derivatives of the
displacement function in terms of integrals of the vector
field~$\mbox{\boldmath$f$}$ along the periodic orbit
$\mbox{\boldmath$\varphi$}_0(t)$~\cite{Gaiko}:
    $$
    d_s(0,\mbox{\boldmath$\mu$}_0)\;=\;\displaystyle\exp\int_{0}^{T_0}\!
    \mbox{\boldmath$\nabla$}\cdot\mbox{\boldmath$f$}
    (\mbox{\boldmath$\varphi$}_0(t),\mbox{\boldmath$\mu$}_0)\:\textrm{d}t-1
    \eqno(3.2)
    \vspace{2mm}
    $$
and
    $$
    d_{\mu_j}(0,\mbox{\boldmath$\mu$}_0)=\\
    $$
    $$
    \frac{-\omega\,_0}
    {\|\mbox{\boldmath$f$}(\mbox{\boldmath$\varphi$}_0(0),
    \mbox{\boldmath$\mu$}_0)\|}\;
    \displaystyle\int_{0}^{T_0}\!\!
    \exp\left(-\!\int_{0}^{t}\!\mbox{\boldmath$\nabla$}\cdot
    \mbox{\boldmath$f$}(\mbox{\boldmath$\varphi$}_0(\tau),
    \mbox{\boldmath$\mu$}_0)\,\textrm{d}\tau\right)
    \mbox{\boldmath$f$}\wedge\mbox{\boldmath$f$}_{\mu_j}
    (\mbox{\boldmath$\varphi$}_0(t),\mbox{\boldmath$\mu$}_0)\:\textrm{d}t\\
    \eqno(3.3)
    $$
for $j=1,\ldots,n,$ where $\omega_0=\pm1$ according to whether $L_0$
is positively or negatively oriented, respectively, and where the wedge
product of two vectors $\mbox{\boldmath$x$}=(x_1,x_2)$ and
$\mbox{\boldmath$y$}=(y_1,y_2)$ in $\textbf{R}^2$ is defined as
    $$
    \mbox{\boldmath$x$}\wedge\mbox{\boldmath$y$}=x_1\,y_2-x_2\,y_1.
    \vspace{-1mm}
    $$
    \par
Similar formulas for $d_{ss}(0,\mbox{\boldmath$\mu$}_0)$ and
$d_{s{\mu_j}}(0,\mbox{\boldmath$\mu$}_0)$ can be derived in terms
of integrals of the vector field $\mbox{\boldmath$f$}$ and its first
and second partial derivatives along $\mbox{\boldmath$\varphi$}_0(t).$
The hypotheses of theorems in the next section will be stated in terms
of conditions on the displacement function $d(s,\mbox{\boldmath$\mu$})$
and its partial derivatives at $(0,\mbox{\boldmath$\mu$}_0)$~\cite{Gaiko}.

\section{Bifurcation surfaces of multiple limit cycles}

In this section, we restate Perko's theorems on the local existence
of $(n\!-\!m\!+\!1)$-dimensional surfaces, $C_m,$ of multiplicity-$m$
limit cycles for the polynomial system (3.1) with
$\mbox{\boldmath$\mu$}\in\textbf{R}^n$ and $n\geq m\geq2.$
These results describe the topological structure of the codimension
$(m\!-\!1)$ bifurcation surfaces $C_m.$ For $m=2,3,4,$
$C_2,$ $C_3,$ and $C_4$ are the familiar fold, cusp, and
swallow-tail bi\-fur\-ca\-tion surfaces; for $m\geq5,$ the
topological structure of the surfaces~$C_m$ is more complex. For
instance, $C_5$ and $C_6$ are the butterfly and wigwam bifurcation
surfaces, respectively~\cite{Perko}. Since the proofs of the theorems
in this section, describing the universal unfolding near a multiple
limit cycles of~(3.1), pa\-ral\-lel the classical proofs of Catastrophe
Theory, we will only state the theorems (see \cite{Perko} for more detail).
    \medskip
    \par
    \textbf{Definition 4.1.}
An $(n\!-\!1)$-dimensional analytic surface
$C_{2}\subset\textbf{R}^n$ is an \emph{$(n\!-\!1)$-dimensional
fold bifurcation surface of multiplicity-two limit cycles of
(3.1) through a point $\mbox{\boldmath$\mu$}_0\in\textbf{R}^n,$}
if for all $\varepsilon>0$ there exists a $\delta>0$ such that
for each $\mbox{\boldmath$\mu$}\in C_{2}$ with
$\|\mbox{\boldmath$\mu$}-\mbox{\boldmath$\mu$}_{0}\|<\delta,$
the system (3.1) has a unique multiplicity-two limit cycle
$L_{\mbox{\boldmath$\mu$}}$ in an $\varepsilon$-neighborhood
of $L_{0}$ and the system (3.1) undergoes a fold bifurcation
at $L_{\mbox{\boldmath$\mu$}};$ i.\,e., for
$\|\mbox{\boldmath$\mu$}-\mbox{\boldmath$\mu$}_{0}\|<\delta,$
$L_{\mbox{\boldmath$\mu$}}$ splits into a simple stable and a
simple unstable limit cycles in an $\varepsilon$-neighborhood of
$L_{0}$ for $\mbox{\boldmath$\mu$}$ on one side of $C_{2}$ and
$L_{\mbox{\boldmath$\mu$}}$ vanishes for $\mbox{\boldmath$\mu$}$
on the other side of $C_{2}.$ \ Cf. Fig.~2.
    \par
\begin{figure}[htb]
\begin{center}
\includegraphics[width=105mm]{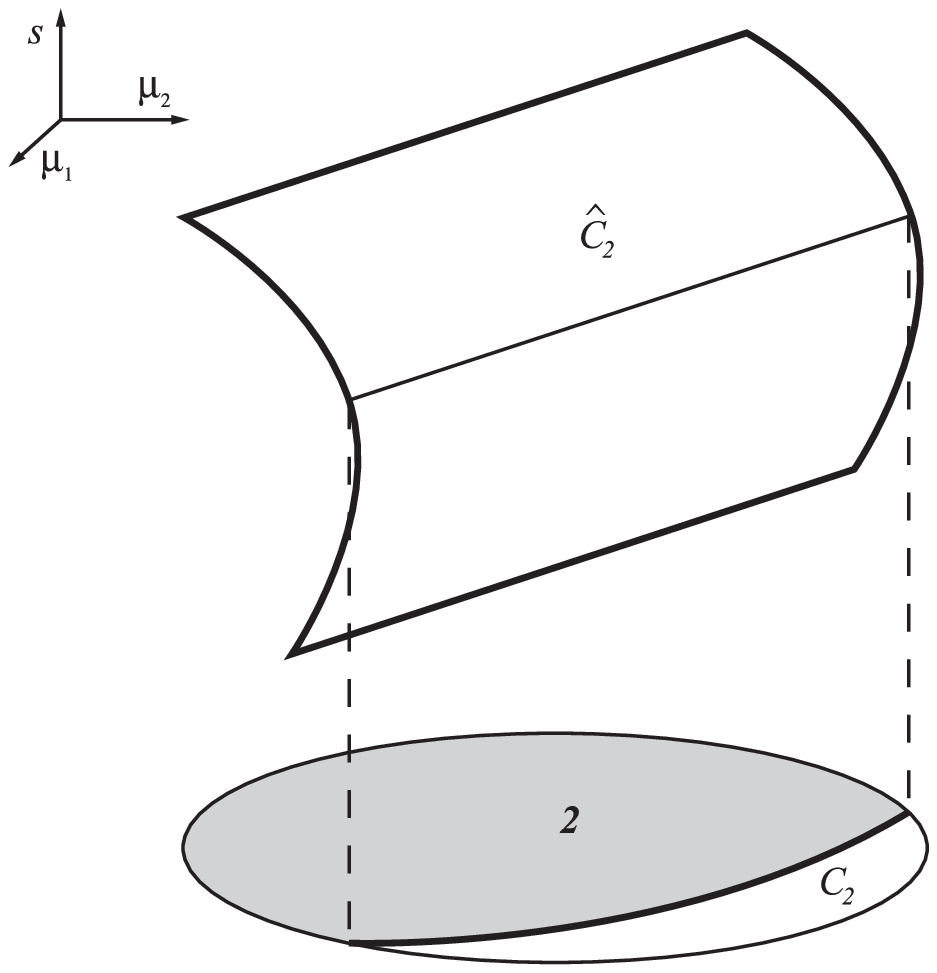}
\vspace{2mm}
    \par
{\small FIG.~2. The fold bifurcation surface.}
\end{center}
\end{figure}
    \par
    \medskip
    \textbf{Theorem 4.1.}
    \emph{Suppose that $n\geq 2,$ that for
    $\mbox{\boldmath$\mu$}=\mbox{\boldmath$\mu$}_0\in\textbf{R}^n$ the
system~(3.1) has a multiplicity-two limit cycle $L_0,$ and that
$d_{\mu_1}(0,\mbox{\boldmath$\mu$}_0)\neq~0.$}
\emph{Then given $\varepsilon>0,$ there is a $\delta>0$
and a unique function $g(\mu_2,\ldots,\mu_n)$ with
$g(\mu_2^{(0)},\ldots,\mu_n^{(0)})=\mu_1^{(0)},$ defined
and analytic for $\vert\mu_2-\mu_2^{(0)}\vert<\delta,$
$\ldots,\vert\mu_n-\mu_n^{(0)}\vert<\delta,$ such that for
$\vert\mu_2-\mu_2^{(0)}\vert<\delta,\ldots,
\vert\mu_n-\mu_n^{(0)}\vert<\delta,$}
\vspace{-1mm}
    $$
    C_2: \quad \mu_1=g(\mu_2,\ldots,\mu_n)\\[-1mm]
    $$
\emph{is an $(n-1)$-dimensional, analytic fold bifurcation surface
of multiplicity-two limit cycles of~(3.1) through the point
$\mbox{\boldmath$\mu$}_0.$}
    \par
    \medskip
    \textbf{Definition 4.2.}
An analytic surface $C_{3}\subset\textbf{R}^n$ is an
\emph{$(n\!-\!2)$-dimen\-si\-onal cusp bifurcation surface of
multiplicity-three limit cycles of (3.1) through a point
$\mbox{\boldmath$\mu$}_0\in\textbf{R}^n,$} if for all
$\varepsilon>0$ there exists a $\delta>0$ such that
for each $\mbox{\boldmath$\mu$}\in C_{3}$ with
$\|\mbox{\boldmath$\mu$}-\mbox{\boldmath$\mu$}_{0}\|<\delta,$
the system (3.1) has a unique multiplicity-three limit cycle
$L_{\mbox{\boldmath$\mu$}}$ in an $\varepsilon$-neighborhood
of $L_{0}$ and the system (3.1) undergoes a cusp bifurcation
at $L_{\mbox{\boldmath$\mu$}};$ i.\,e., $C_{3}$ is the intersection
of two $(n\!-\!1)$-dimensional fold bifurcation surfaces of
multiplicity-two limit cycles of (3.1), $C_{2}^{\pm},$ which
intersect in a cusp along $C_{3};$ for
$\|\mbox{\boldmath$\mu$}-\mbox{\boldmath$\mu$}_{0}\|<\delta$
and for $\mbox{\boldmath$\mu$}$ in the cuspidal region between
$C_{2}^{+}$ and $C_{2}^{-}$ (shaded in Fig.~3), the system (3.1)
has three simple limit cycles in an $\varepsilon$-neighborhood of $L_{0};$
and for $\|\mbox{\boldmath$\mu$}-\mbox{\boldmath$\mu$}_{0}\|<\delta$ and
$\mbox{\boldmath$\mu$}$ outside the cuspidal region, the system (3.1) has
one simple limit cycle in an $\varepsilon$-neighborhood of~$L_{0}.$ Cf.~Fig.~3.
    \par
\begin{figure}[htb]
\begin{center}
\includegraphics[width=105mm]{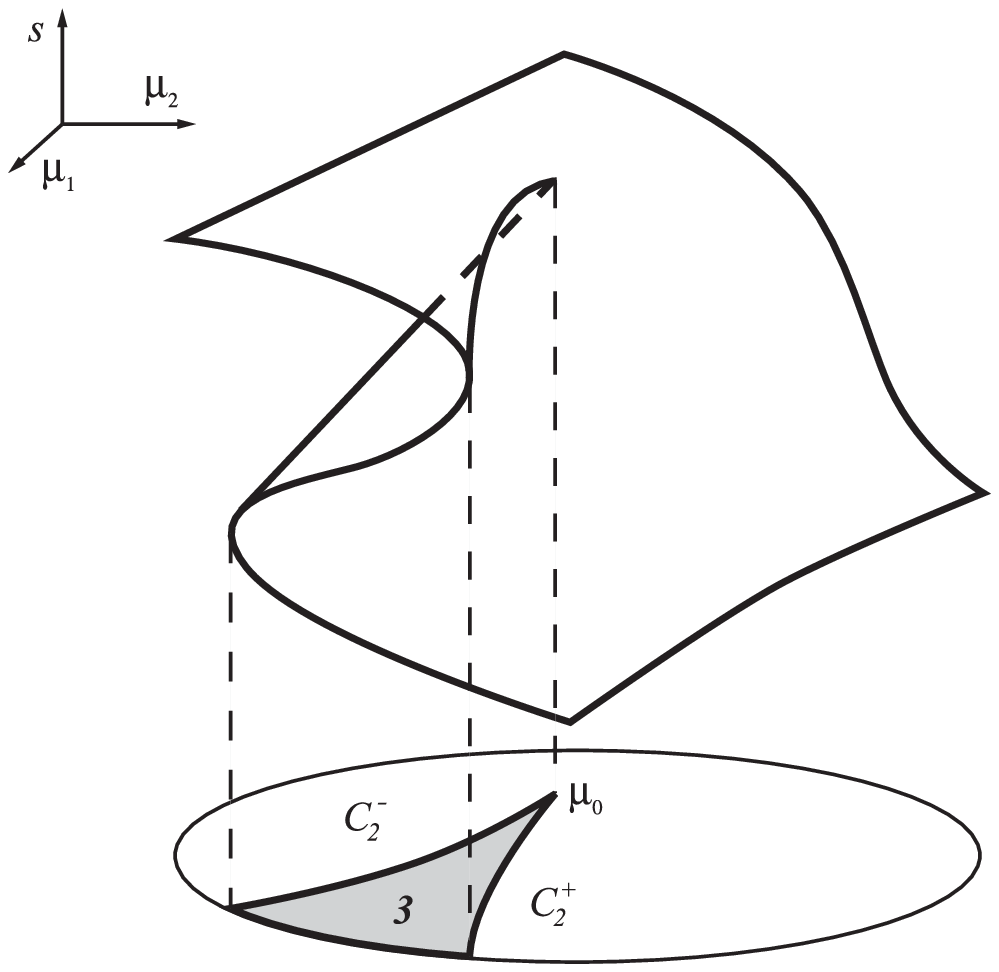}
    \vspace{4mm}
    \par
{\small FIG.~3. The cusp bifurcation surface.}
\end{center}
\end{figure}
    \par
    \medskip
    \textbf{Theorem 4.2.}
    \emph{Suppose that $n\geq3,$ that for $\mbox{\boldmath$\mu$}=
\mbox{\boldmath$\mu$}_0\in\textbf{R}^n$ the system (3.1) has
a multiplicity-three limit cycle $L_0,$ that
$d_{\mu_1}(0,\mbox{\boldmath$\mu$}_0)\neq0,$
$d_{r{\mu_1}}(0,\mbox{\boldmath$\mu$}_0)\neq0$
and for $j=2,\ldots,n,$}
    $$
    \Delta_j\equiv\frac{\partial(d, d_r)}
    {\partial(\mu_1,\mu_j)}(0,\mbox{\boldmath$\mu$}_0)\neq0.
    $$
    \par
    \emph{Then given $\varepsilon>0,$ there is a $\delta>0$ and
constants $\omega_j=\pm1$ for $j=2,\ldots,n,$ and there exist
unique functions $h_1(\mu_2,\ldots,\mu_n),$
$h_2(\mu_2,\ldots,\mu_n)$ and $g^{\pm}(\mu_2,\ldots,\mu_n)$
with $h_1(\mu_2^{(0)},\ldots,\mu_n^{(0)})=\mu_1^{(0)},$
$h_2(\mu_2^{(0)},\ldots,\mu_n^{(0)})=\mu_1^{(0)}$ and
$g^{\pm}(\mu_2^{(0)},\ldots,\mu_n^{(0)})=\mu_1^{(0)},$
where $h_1$ and $h_2$ are defined and analytic for
$\vert\mu_j-\mu_j^{(0)}\vert<\delta,$ $j=2,\ldots,n,$
and $g^{\pm}$ are defined and continuous for
$0\leq\sigma_j(\mu_j-\mu_j^{(0)})<\delta$ and analytic
for $0<\omega_j(\mu_j-\mu_j^{(0)})<\delta,$ $j=2,\ldots,n$
such that}
    $$
    C_3: \quad
    \left
    \{
    \begin{array}{rl}
    \mu_{1}=h_{1}(\mu_{2},\ldots,\mu_{n})\\
    \mu_{1}=h_{2}(\mu_{2},\ldots,\mu_{n})\\
    \end{array}
    \right.
    $$
\emph{is an $(n-2)$-dimensional, analytic, cusp bifurcation
surface of multi\-pli\-city-three limit cycles of~(3.1)
through the point $\mbox{\boldmath$\mu$}_0$ and}
    $$
    C_2^{\pm}: \quad
    \mu_{1}=g^{\pm}(\mu_{2},\ldots,\mu_{n})
    $$
\emph{are two $(n-1)$-dimensional, analytic, fold bifurcation
surfaces of multi\-pli\-city-two limit cycles of (3.1) which
intersect in a cusp along
$C_3.$}
    \medskip
    \par
    \textbf{Definition 4.3.}
An analytic surface $C_{4}\subset\textbf{R}^n$ is an
\emph{$(n\!-\!3)$-dimen\-si\-onal swallow-tail bifurcation surface
of multiplicity-four limit cycles of (3.1) through a point
$\mbox{\boldmath$\mu$}_0\in\textbf{R}^n,$} if for all
$\varepsilon>0$ there exists a $\delta>0$ such that for
each $\mbox{\boldmath$\mu$}\in C_{4}$ with
$\|\mbox{\boldmath$\mu$}-\mbox{\boldmath$\mu$}_{0}\|<\delta,$
the system (3.1) has a unique multiplicity-four limit cycle
$L_{\mbox{\boldmath$\mu$}}$ in an $\varepsilon$-neighborhood
of $L_{0}$ and the system (3.1) undergoes a swallow-tail
bifurcation at $L_{\mbox{\boldmath$\mu$}};$ i.\,e., $C_{4}$ is
the intersection of two $(n\!-\!2)$-dimensional cusp bifurcation
surfaces of multiplicity-three limit cycles $C_{3}^{\pm}$ which
intersect in a cusp along $C_{4};$ furthermore, there are three
$(n\!-\!1)$-dimensional fold bifurcation surfaces of
multiplicity-two limit cycles of (3.1),
$C_{2}^{(i)},$ $i=0,1,2,$ such that $C_{2}^{(0)}$ and
$C_{2}^{(1)}$ intersect in a cusp along $C_{3}^{+},$ $C_{2}^{(0)}$
and $C_{2}^{(2)}$ intersect in a cusp along $C_{3}^{-},$ and
$C_{2}^{(1)}$ and $C_{2}^{(2)}$ intersect along an
$(n\!-\!2)$-dimensional surface on which (3.1) has
two multiplicity-two limit cycles; finally, for
$\|\mbox{\boldmath$\mu$}-\mbox{\boldmath$\mu$}_{0}\|<\delta$
and for $\mbox{\boldmath$\mu$}$ in the swallow-tail region
(shaded in Fig.~4), the system (3.1) has four simple limit cycles
in an $\varepsilon$-neighborhood of $L_{0};$ for
$\|\mbox{\boldmath$\mu$}-\mbox{\boldmath$\mu$}_{0}\|<\delta$
and $\mbox{\boldmath$\mu$}$ above the surfaces $C_{2}^{(i)},$
$i=0,1,2,$ the system (3.1) has two simple limit cycles in an
$\varepsilon$-neighborhood of $L_{0};$ and for
$\|\mbox{\boldmath$\mu$}\!-\!\mbox{\boldmath$\mu$}_{0}\|\!<\!\delta$
and $\mbox{\boldmath$\mu$}$ below the surfaces $C_{2}^{(i)},$
$i=0,1,2,$ the system (3.1) has no limit cycles in an
$\varepsilon$-neighborhood of $L_{0}.$ Cf.~Fig.~4.
    \par
\begin{figure}[htb]
\begin{center}
\includegraphics[width=105mm]{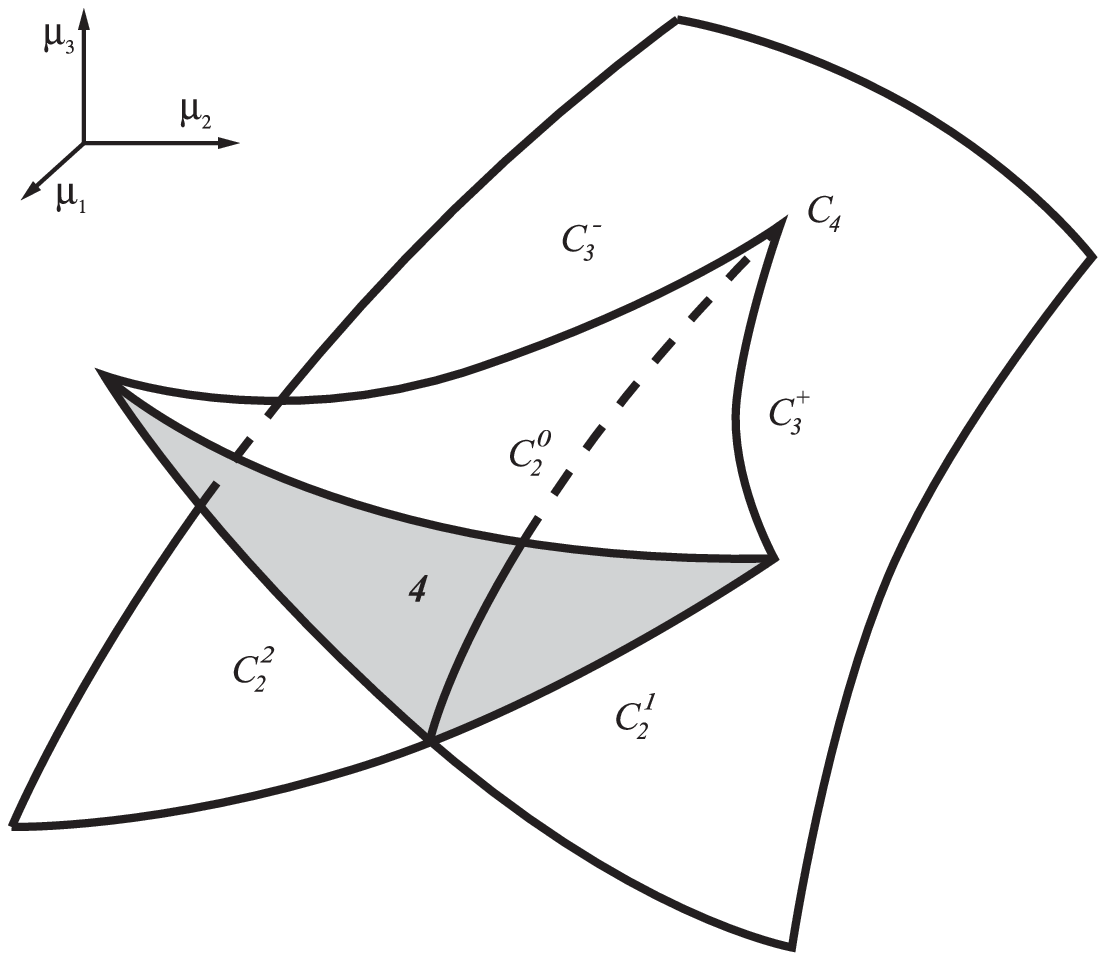}
    \par
{\small FIG.~4. The swallow-tail bifurcation surface.}
\end{center}
\end{figure}
    \par
    \bigskip
    \textbf{Theorem 4.3.}
    \emph{Suppose that $n\geq4,$ that for
$\mbox{\boldmath$\mu$}=\mbox{\boldmath$\mu$}_0\in\textbf{R}^n$
the system (3.1) has a multiplicity-four limit cycle $L_0,$
that $d_{\mu_1}(0,\mbox{\boldmath$\mu$}_0)\neq0,$
$d_{r{\mu_1}}(0,\mbox{\boldmath$\mu$}_0)\neq0,$
$d_{rr{\mu_1}}(0,\mbox{\boldmath$\mu$}_0)\neq0,$
and that for $j=2,\ldots,n,$}
    $$
    \frac{\partial(d,d_r)}{\partial(\mu_1,\mu_j)}
    (0,\mbox{\boldmath$\mu$}_0)\neq 0, \quad
    \frac{\partial(d,d_{rr})}{\partial(\mu_1,\mu_j)}
    (0,\mbox{\boldmath$\mu$}_0)\neq 0, \quad
    \frac{\partial(d_r,d_{rr})}{\partial(\mu_1,\mu_j)}
    (0,\mbox{\boldmath$\mu$}_0)\neq0.\\[-2mm]
    $$
    \par
    \emph{Then given $\varepsilon>0,$ there is a $\delta>0$ and
constants $\omega_{jk}=\pm1$ for $j=2,\ldots,n,$ $k=1,2,$ and
there exist unique functions $g_i(\mu_2,\ldots,\mu_n),$
$h_k^{\pm}(\mu_2,\ldots,\mu_n)$ and $F_i(\mu_2,\ldots,\mu_n),$
with $g_i(\mu_2^{(0)},\ldots,\mu_n^{(0)})=
h_k^{\pm}(\mu_2^{(0)},\ldots,\mu_n^{(0)})=
F_i(\mu_2^{(0)},\ldots,\mu_n^{(0)})\!=\!\mu_1^{(0)}\!,$ for
$i=0,1,2$ and $k=1,2,$ where $F_i$ is defined and analytic for
$i=0,1,2,$ and $\vert\mu_j-\mu_j^{(0)}\vert<\delta,$
$j=2,\ldots,n,$ $h_k^{\pm}$ are defined and continuous for
$0\leq\omega_{jk}(\mu_j-\mu_j^{(0)})<\delta$ and analytic for
$0<\omega_{jk}(\mu_j-\mu_j^{(0)})<\delta,$ $j=2,\ldots,n,$
$k=1,2,$ and for $i=0,1,2,$ $g_i$ is defined and analytic in the
cuspidal region between the surfaces
$\mu_1=h_2^{\pm}(\mu_2,\ldots,\mu_n),$ which intersect in a cusp,
and $g_i$ is continuous in the closure of that region, such that}
    $$
    C_4: \quad
    \left
    \{
    \begin{array}{rl}
    \mu_{1}=F_{0}(\mu_{2},\ldots,\mu_{n})\\
    \mu_{1}=F_{1}(\mu_{2},\ldots,\mu_{n})\\
    \mu_{1}=F_{2}(\mu_{2},\ldots,\mu_{n})\\
    \end{array}
    \right.
    $$
\emph{is an $(n-3)$-dimensional, analytic, swallow-tail bifurcation
surface of multi\-pli\-city-four limit cycles of (3.1) through the
point $\mbox{\boldmath$\mu$}_0$ which is the intersection of two
$(n-2)$-dimensional, analytic, cusp bifurcation surfaces of
multiplicity-three limit cycles of (3.1),}
    $$
    C_3^{\pm}: \quad
    \left
    \{
    \begin{array}{rl}
    \mu_{1}=h_{1}^{\pm}(\mu_{2},\ldots,\mu_{n})\\
    \mu_{1}=h_{2}^{\pm}(\mu_{2},\ldots,\mu_{n})\\
    \end{array}
    \right.
    $$
\emph{which intersect in a cusp along $C_4;$ furthermore,
$C_3^{+}=C_2^{(0)} \bigcap C_2^{(1)}$ and $C_3^{-}=C_2^{(0)}
\bigcap C_2^{(2)}$ where for $i=0,1,2,$}
    $$
    C_2^{i}: \quad
    \mu_{1}=g_{i}(\mu_{2},\ldots,\mu_{n})
    $$
\emph{are $(n-1)$-dimensional, analytic, fold bifurcation surfaces
of multiplicity-two limit cycles of (3.1)  which intersect in cusps
along $C_3^{\pm}$ and in an $(n-2)$-dimensional, analytic surface
$C_2^{(1)} \bigcap C_2^{(2)}$ on which (3.1) has two multiplicity-two
limit cycles} (Fig.~4 and Fig.~5).
    \par
\begin{figure}[htb]
\begin{center}
\includegraphics[width=100mm]{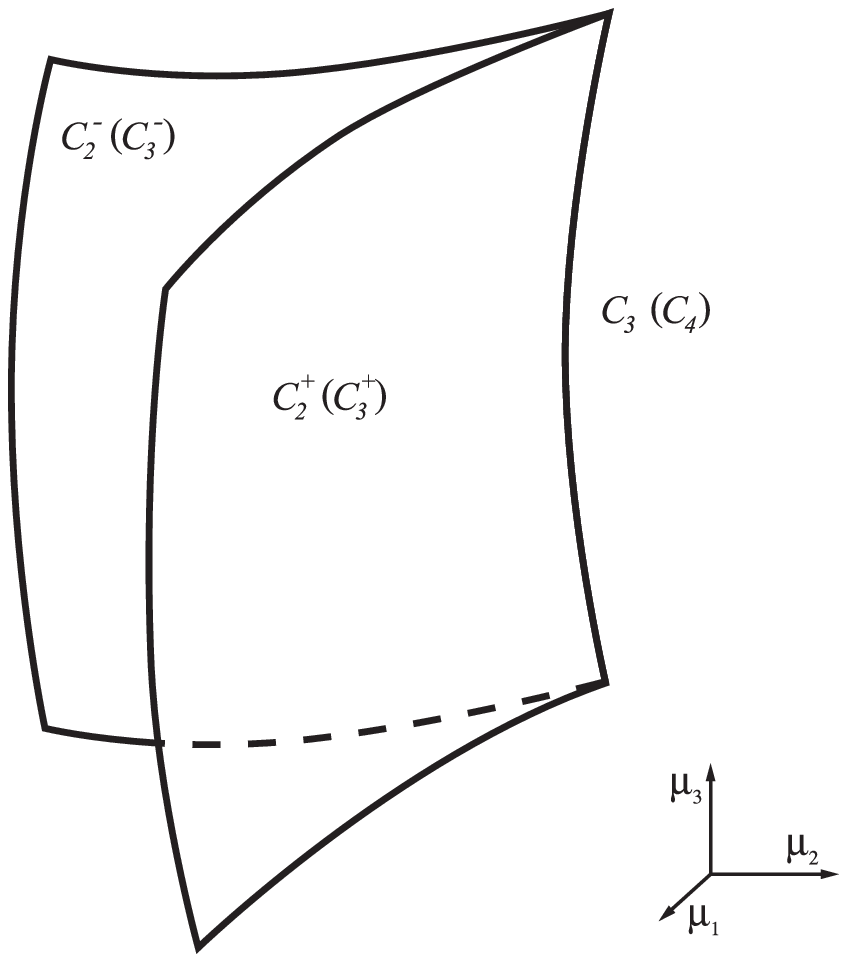}
    \par
{\small FIG.~5. The bifurcation curve (one-parameter family)
of multiple limit cycles.}
\end{center}
\end{figure}
    \medskip
    \par
Based on Theorems~3.1,~3.2, the following ge\-ne\-ra\-li\-za\-tion of
Theorems~4.1\,--\,4.3 can be proved on induction~\cite{Perko}.
    \par
    \medskip
    \textbf{Theorem 4.4.}
    \emph{Given $m\geq2.$ Suppose that $n\geq m,$ that for
$\mbox{\boldmath$\mu$}=\mbox{\boldmath$\mu$}_0\in\textbf{R}^n$ the
polynomial system~(3.1) has a multiplicity-$m$ limit cycle $L_0,$ that}
    $$
    \frac{\partial d}{\partial\mu_1}\,(0,\mbox{\boldmath$\mu$}_0)\neq0,
        \quad
    \frac{\partial d_r}{\partial\mu_1}\,(0,\mbox{\boldmath$\mu$}_0)\neq0,\;
        \ldots,
        \quad
    \frac{\partial d_r^{(m-2)}}{\partial\mu_1}\,(0,\mbox{\boldmath$\mu$}_0)\neq0,
    $$
\emph{and that}
    $$
    \frac{\partial(d_r^{(i)},d_r^{(j)})}{\partial(\mu_1,\mu_k)}\,
    (0,\mbox{\boldmath$\mu$}_0)\neq0
    \vspace{4mm}
    $$
\emph{for $i,j=0,\ldots,m-2$ with $i\neq j$ and $k=2,\ldots,n.$}
    \par
    \emph{Then given $\varepsilon>0$ there is a $\delta>0$ such that
for $\|\mbox{\boldmath$\mu$}-\mbox{\boldmath$\mu$}_0\|<\delta,$
the system (3.1) has}
    \par
$(1)$ \ \emph{a unique $(n-m+1)$-dimensional analytic surface $C_m$
of multi\-pli\-city-$m$ limit cycles of (3.1) through the point
$\mbox{\boldmath$\mu$}_0;$}
    \par
$(2)$ \ \emph{two $(n-m+2)$-dimensional analytic surfaces $C_{m-1}$
of multiplicity-$(m\!-\!1)$ limit cycles of (3.1) through the point
$\mbox{\boldmath$\mu$}_0$ which intersect in a cusp along $C_m;$}
    \par
    \vspace{-4mm}
    \dots
    \vspace{-2mm}
    \par
$(j)$ \ \emph{exactly $j,$ $(n-m+j)$-dimensional analytic surfaces
$C_{m-j+1}$ of multiplicity-$(m-j+1)$ limit cycles of (3.1) through
the point $\mbox{\boldmath$\mu$}_0$ which intersect pairwise in cusps
along the bifurcation surfaces $C_{m-j+2};$}
    \par
    \vspace{-4mm}
    $\dots$
    \vspace{-2mm}
    \par
$(m-1)$ \ \emph{exactly $(m-1),$ $(n-1)$-dimensional analytic fold
bifurcation surfaces $C_2$ of multiplicity-two limit cycles of (3.1)
through the point $\mbox{\boldmath$\mu$}_0$ which intersect pairwise
in a cusp along the $(n-2)$-dimensional cusp bifurcation surfaces $C_3.$}

\section{Singular points of the FitzHugh--Nagumo system}

The study of singular point of the system (1.3) will use two index
theorems by H.\,Poincar\'{e}, see \cite{BL}. But first let us define
the Poincar\'{e} index~\cite{BL}.
    \medskip
    \par
    \textbf{Definition 5.1.}
Let $S$ be a simple closed curve in the phase plane not passing
through a singular point of the system
    $$
    \dot{x}=P(x,y), \quad \dot{y}=Q(x,y),
    \eqno(5.1)
    $$
where $P(x,y)$ and $Q(x,y)$ are continuous functions (for example,
polynomials), and $M$ be some point on $S.$ If the point $M$ goes
around the curve $S$ in positive direction (counterclockwise) one
time, then the vector coinciding with the direction of a tangent
to the trajectory passing through the point $M$ is rotated through
the angle $2\pi j$ $(j=0,\pm1,\pm2,\ldots).$ The integer $j$ is
called the \emph{Poincar\'{e} index} of the closed curve $S$
relative to the vector field of system~(5.1) and has the
expression
    $$
    j=\frac{1}{2\pi}\oint_S\frac{P~dQ-Q~dP}{P^2+Q^2}.\\[-2mm]
    $$
    \par
According to this definition, the index of a node or a focus, or a
center is equal to $+1$ and the index of a saddle is $-1.$
    \par
    \medskip
    \textbf{Theorem 5.1 (First Poincar\'{e} Index Theorem).}
    \emph{If $N,$ $N_f,$ $N_c,$ and $C$ are respectively the number
of nodes, foci, centers, and saddles in a finite part of the phase
plane and $N'$ and $C'$ are the number of nodes and saddles at
infinity, then it is valid the formula}
    $$
    N+N_f+N_c+N'=C+C'+1.
    $$
    \par
    \textbf{Theorem 5.2 (Second Poincar\'{e} Index Theorem).}
    \emph{If all singular points are simple, then along an isocline
without multiple points lying in a Poincar\'{e} hemisphere which is
obtained by a stereographic projection of the phase plane, the
singular points are distributed so that a saddle is followed by a
node or a focus, or a center and vice versa. If two points are
separated by the equator of the Poincar\'{e} sphere, then a saddle
will be followed by a saddle again and a node or a focus, or a
center will be followed by a node or a focus, or a center.}
    \medskip
    \par
Consider the system (1.3). Its finite singularities are determined by the algebraic system
   $$
    \begin{array}{l}
(\gamma\,\delta-1)\,y+(\gamma-a)\,x+b\,x^2-c\,x^3=0,
    \\
x-\delta\,y=0\,.
    \end{array}
    \eqno(5.2)
    $$
From (5.2) \cite{BL}, \cite{Gaiko}, \cite{rzh}, we will get a singular point $(0,0)$ of antisaddle-type
(i.\,e., a node, a focus, or a center) and at most two points (a saddle and an antisaddle or,
if one, a saddle-node) defined by the condition
    $$
\displaystyle c\,x^2-b\,x-\frac{\gamma\,\delta-1}{\delta}-\gamma+a=0,
    \quad
y=\frac{x}{\delta}\,.
    \eqno(5.3)
    $$
To get singular points at infinity, consider the corresponding differential equation
    $$
\frac{dy}{dx}=\frac{x-\delta\,y}{(\gamma\,\delta-1)\,y+(\gamma-a)\,x+b\,x^2-c\,x^3}\,.\\[2mm]
    \eqno(5.4)
    $$
Dividing the numerator and denominator of the right-hand side of (5.4) first by $x^{3}$
$(x\neq0),$ denoting $y/x$ by $u,$ and then by $y^{3}$ $(y\neq0),$ denoting $x/y$ by $v,$
we will get two infinite singularities: $u=0$ (a simple node in the direction of the $x$-axis)
and $v^{3}=0$ (a triple saddle in the direction of the $y$-axis), see~\cite{BL},~\cite{Gaiko}.

\section{Global limit cycle bifurcations in the system}

To investigate global limit cycle bifurcations in the system (1.3), we will use the results of
the previous sections and will apply the method developed in \cite{bog}, \cite{bg}, \cite{Gaiko}--\cite{gai6}.
The sense of this method is to obtain the simplest (well-known) system by vanishing some parameters (usually
field rotation parameters) of the original system and then to input these parameters successively
one by one studying the dynamics of limit cycles on the whole phase plane.
    \par
Let us study rotation properties of the parameters of (1.3). Applying the definition of a field
rotation parameter (i.\,e., a parameter which rotates the field in one direction \cite{BL},
\cite{Gaiko}, \cite{Perko}), to the system (1.3) written in the form
    $$
\dot{x}=R(x,y)+\gamma\,Q(x,y)\equiv P(x,y), \quad \dot{y}=Q(x,y),
    \eqno(6.1)
    $$
where $R(x,y)=-y-a\,x+b\,x^2-c\,x^3$ and $Q(x,y)=x-\delta\,y,$ calculate the corresponding
determinant for the parameter $\gamma\!:$
    $$
\Delta_{\gamma}=PQ'_{\gamma}-QP'_{\gamma}=-Q^{2}\leq0.
    \eqno(6.2)
    $$
It follows from (6.2) that on increasing $\gamma$ the vector field of~(1.3)
is rotated in negative direction (clockwise) on the whole phase plane of~(1.3).
    \par
For $\delta=0,$ we will have a system
    \vspace{-1mm}
    $$
    \begin{array}{l}
\dot{x}=-y+(\gamma-a)\,x+b\,x^2-c\,x^3\equiv P(x,y),
    \\
\dot{y}=x\equiv Q(x,y).
    \end{array}
    \eqno(6.3)
    $$
Applying the definition of a field rotation parameter to this system for the parameters $a$ and $c,$
respectively, we will get the following determinants:
    $$
\Delta_{a}=PQ'_{a}-QP'_{a}=x^{2}\geq0,\\
    \eqno(6.4)
    $$
    $$
\Delta_{c}=PQ'_{c}-QP'_{c}=x^{4}\geq0.\\[2mm]
    \eqno(6.5)
    $$
It follows from (6.4) and (6.5) that on increasing $a$ or $c$ the vector field of~(6.3)
is rotated in positive direction (counterclockwise) on the whole phase plane of~(6.3).
    \par
For the study of multiple limit cycle bifurcations, we will use also two theorems by
L.\,Perko (see~\cite{Perko}) which are formulated for the polynomial system (3.1).
    \par
    \medskip
\noindent\textbf{Theorem 6.1 (Wintner--Perko termination principle).}
    \emph{Any one-para\-me\-ter fa\-mi\-ly of multip\-li\-city-$m$ limit cycles
of the relatively prime polynomial system (3.1) can be extended in a unique way
to a maximal one-parameter family of multiplicity-$m$ limit cycles of (3.1)
which is either open or cyclic.}
    \par
\emph{If it is open, then it terminates either as the parameter or the limit cycles
become unbounded; or, the family terminates either at a singular point of (3.1),
which is typically a fine focus of multiplicity~$m,$ or on a (compound)
separatrix cycle of (3.1), which is also typically of multiplicity~$m.$}
    \medskip
    \par
The proof of this principle for the general polynomial system (3.1)
with a vector parameter $\mbox{\boldmath$\mu$}\in\textbf{R}^n$
parallels the proof of the pla\-nar termination principle for
the system
    $$
    \vspace{1mm}
    \dot{x}=P(x,y,\lambda),
        \quad
    \dot{y}=Q(x,y,\lambda)\\[2mm]
    \eqno(6.6)
    $$
with a scalar parameter $\lambda\in\textbf{R},$ since there is no loss
of generality in assuming that system (3.1) is parameterized by a scalar
parameter $\lambda$ (see \cite{Gaiko}, \cite{Perko}).
    \par
In particular, if $\lambda$ is a field rotation parameter of (6.6),
the following Perko's theorem on monotonic families of multiple
limit cycles is valid~\cite{Perko}.
    \par
    \medskip
\noindent\textbf{Theorem 6.2.}
    \emph{If $L_{0}$ is a nonsingular multiple limit cycle of (6.6) for
$\lambda=\lambda\,_{0},$ then  $L_{0}$ belongs to a one-parameter family
of limit cycles of (6.6); furthermore:}
    \par
1)~\emph{if the multiplicity of $L_{0}$ is odd, then the family
either expands or contracts mo\-no\-to\-ni\-cal\-ly as $\lambda$
increases through $\lambda_{0};$}
    \par
2)~\emph{if the multiplicity of $L_{0}$ is even, then $L_{0}$
bi\-fur\-cates into a stable and an unstable limit cycle as
$\lambda$ varies from $\lambda_{0}$ in one sense and $L_{0}$
dis\-ap\-pears as $\lambda$ varies from $\lambda_{0}$ in the
opposite sense; i.\,e., there is a fold bifurcation at
$\lambda_{0}.$}
    \par
Using these theorems and the results of the previous sections
and applying the field rotation parameters of the systems~(1.3)
and~(6.3), we will prove the following theorem.
	\par
    \medskip
\noindent \textbf{Theorem 6.3.}
\emph{The FitzHugh--Nagumo system (1.3) can have at most two limit cycles.}
	\par
    \medskip
\noindent\textbf{Proof.} First let us prove that system (1.3) can have at least two limit cycles
supposing that all of the parameters of (1.3) are nonnegative (in the natural sense of the model).
All other cases can be considered in a similar way.
    \par
Let the parameters $a,$ $c,$ $\gamma,$ $\delta$ of (1.3) vanish and consider the quadratic system
    $$
\dot{x}=-y+b\,x^2, \quad \dot{y}=x.
    \eqno(6.7)
    $$
This is a reversible system. It has a center at the origan $O$ and cannot have limit cycles \cite{Gaiko}.
We will input the parameters $\gamma,$ $a,$ $c,$ and $\delta$ successively one by one into the system (6.7).
    \par
Inputting a positive parameter $\gamma,$ we will get a~system
    $$
\dot{x}=-y+\gamma\,x+b\,x^2, \quad \dot{y}=x,
    \eqno(6.8)
    $$
the vector field of which is rotated in negative direction (clockwise) on the whole phase plane of~(6.8).
The origin of~(6.8) becomes an unstable focus (or a node).
    \par
Inputting a positive parameter $a$ into~(6.8), the vector field of the system
    $$
\dot{x}=-y+(\gamma-a)\,x+b\,x^2, \quad \dot{y}=x
    \eqno(6.9)
    $$
will be rotated in positive direction (counterclockwise). For $a=\gamma,$ the origin becomes weak and changes
the character of stability on further increasing~$a.$ The Andronov--Hopf bifurcation occurs for $a=\gamma,$
and an unstable limit cycle, $\Gamma_{1},$ will appear from the origin \cite{BL}, \cite{Gaiko}.
    \par
Inputting a positive parameter $c$ into (6.9), we will get a cubic system
    $$
\dot{x}=-y+(\gamma-a)\,x+b\,x^2-c\,x^3, \quad \dot{y}=x,
    \eqno(6.9)
    $$
the vector field of which is also rotated in positive direction on the whole phase plane of~(6.9).
The structure and the character of stability of infinite singularities will be changed, and a stable limit,
$\Gamma_{2},$ surrounding $\Gamma_{1}$ will appear immediately from infinity in this case. On further increasing
the parameter~$c,$ the limit cycles $\Gamma_{1}$ and $\Gamma_{2}$ combine a semi-stable limit, $\Gamma_{12},$
which then disappears in a ``trajectory concentration'' \cite{BL}, \cite{Gaiko}.
    \par
If to input a positive parameter $\delta$ into (6.9), we will have again the original system (1.3).
On further increasing $\delta,$ a saddle-node appears in the first quadrant of the phase plane. It splits then
in two singular points: a saddle $S$ and an antisaddle $A.$ Without loss of generality, we can fix the parameter
$\delta,$ fixing the positions of the finite singularities $O,$ $S,$ $A,$ and consider the system (1.3)
with a positive parameter~$\gamma$ which rotates the vector field of~(1.3) on the whole phase plane.
    \par
So, consider the original system (1.3) with a positive parameter~$\gamma.$ On increasing this parameter, the stable
nodes $O$ and $A$ becomes first stable foci, then they change the character of their stability, becoming unstable foci.
At these Andronov--Hopf bifurcations \cite{BL}, \cite{Gaiko}, stable limit cycles will appear from the foci $O$ and
$A.$ On further increasing~$\gamma,$ the limit cycles will expand and will disappear in small separatrix loops of
the saddle $S.$ If these loops are formed simultaneously, we will have a so-called eight-loop separatrix cycle.
In this case, a big stable limit surrounding three singular points, $O,$ $S,$ and $A,$ will appear from the
eight-loop separatrix cycle after its destruction, expanding to infinity on increasing~$\gamma.$
If a small loop is formed earlier, for example, around the point $O$ $(A),$ then, on increasing~$\gamma,$
a big loop formed by two lower (upper) adjoining separatrices of the saddle~$S$ and surrounding the points $O$
and $A$ will appear. After its destruction, we will have simultaneously a big limit cycle surrounding
three singular points, $O,$ $S,$ $A,$ and a small limit cycle surrounding the point $A$ $(O).$ Thus,
we have proved that system (1.3) can have at least two limit cycles, see also \cite{rzh} for more detail.
    \par
Let us prove now that this system has at most two limit cycles. The proof is carried out by contradiction applying
Catastrophe Theory, see \cite{Gaiko}, \cite{Perko}. Consider the system (1.3) with three field rotation parameters:
$a,$ $c,$ and $\gamma$ (the parameters $b$ and $\delta$ can be fixed, since they do not generate limit cycles).
Suppose that (1.3) has three limit cycles surrounding the only point, $O,$ at~the origin. Then we get into some
domain of the parameters $a$ $c,$ and $\gamma$ being restricted by definite con\-di\-tions on two other parameters,
$b$ and $\delta.$ This domain is bounded by two fold bifurcation surfaces forming a cusp bifurcation surface of
multiplicity-three limit cycles in the space of the pa\-ra\-me\-ters $a,$ $c,$ and $\gamma$ \cite{Gaiko},~\cite{Perko}.
    \par
The cor\-res\-pon\-ding maximal one-parameter family of multiplicity-three limit cycles cannot be cyclic, otherwise there
will be at least one point cor\-res\-pon\-ding to the limit cycle of multi\-pli\-ci\-ty four (or even higher) in the
parameter space. Extending the bifurcation curve of multi\-pli\-ci\-ty-four limit cycles through this point and parameterizing
the corresponding maximal one-parameter family of multi\-pli\-ci\-ty-four limit cycles by a field rotation para\-me\-ter,
$\gamma,$ according to Theorem~6.2, we will obtain two monotonic curves of, respectively, multi\-pli\-ci\-ty-three and
one limit cycles which, by the Wintner--Perko termination principle (Theorem~6.1), terminate either at the point $O$ or
on an infinite separatrix cycle surrounding this point. Since we know at least the cyclicity of the singular point which
is equal to two (see \cite{rzh}), we have got a contradiction with the termination principle stating that the multiplicity
of limit cycles cannot be higher than the multi\-pli\-ci\-ty (cyclicity) of the singular point in which they terminate.
    \par
If the maximal one-parameter family of multiplicity-four limit cycles is not cyclic, using the same principle (Theorem~6.1),
this again contradicts the cyclicity of the origin (see \cite{rzh}) not admitting the multiplicity of limit cycles to be higher
than two. This contradiction completes the proof in the case of one singular point on the phase plane.
    \par
Suppose that the system (1.3) has three finite singularities, $O,$ $S,$ $A,$ and two small limit cycles around,
e.\,g., the point $O$ (the case when the limit cycles surround the point $A$ is considered in a similar way).
Then we get into some domain in the space of the parameters $a,$ $c,$ and $\gamma$ which is bounded by
a fold bifurcation surface of multiplicity-two limit cycles \cite{Gaiko},~\cite{Perko}.
    \par
The cor\-res\-pon\-ding maximal one-parameter family of multiplicity-two limit cycles cannot be cyclic, otherwise there
will be at least one point cor\-res\-pon\-ding to the limit cycle of multi\-pli\-ci\-ty three (or even higher) in the
parameter space. Extending the bifurcation curve of multi\-pli\-ci\-ty-three limit cycles through this point and
parameterizing the corresponding maximal one-parameter family of multi\-pli\-ci\-ty-three limit cycles by a field
rotation para\-me\-ter, $\gamma,$ according to Theorem~6.2, we will obtain a monotonic curve which, by the Wintner--Perko
termination principle (Theorem~6.1), terminates either at the point $O$ or on some separatrix cycle surrounding this point.
Since we know at least the cyclicity of the singular point which is equal to one in this case \cite{rzh},
we have got a contradiction with the termination principle (Theorem~6.1).
    \par
If the maximal one-parameter family of multiplicity-two limit cycles is not cyclic, using the same principle (Theorem~6.1),
this again contradicts the cyclicity of $O$ (see \cite{rzh}) not admitting the multiplicity of limit cycles higher than one.
Moreover, it also follows from the termination principle that either an ordinary (small) separatrix loop or a big loop, or
an eight-loop cannot have the multiplicity (cyclicity) higher than one in this case. Therefore, according to the same principle,
there are no more than one limit cycle in the exterior domain surrounding all three finite singularities, $O,$ $S,$ and $A.$
    \par
Thus, taking into account all other possibilities for limit cycle bifurcations (see \cite{rzh}), we conclude that system~(1.3)
cannot have either a multiplicity-three limit cycle or more than two limit cycles in any configuration. The theorem is proved.
\qquad $\Box$

\end{document}